\documentclass{amsart}

\usepackage[a4paper,hmargin=3cm,vmargin=4cm]{geometry}
\usepackage{amsfonts,amssymb,amscd,amstext}
\usepackage{hyperref}
\usepackage{mathpazo}
\usepackage{graphicx}

\let\geq=\geqslant
\let\leq=\leqslant
\let\ptl=\partial
\let\Sg=\Sigma
\let\sg=\sigma

\let\Om=\Omega

\newcommand{\rr}{\mathbb{R}}
\newcommand{\rrn}{\mathbb{R}^{2n+1}}
\newcommand{\hh}{\mathbb{H}}
\newcommand{\hhn}{\mathbb{H}^n}
\newcommand{\pp}{\mathcal{P}}
\newcommand{\nuh}{\nu_{H}}
\newcommand{\var}{\varphi}
\newcommand{\sub}{\subset}
\newcommand{\subeq}{\subseteq}
\newcommand{\escpr}[1]{\left< #1\right>}

\DeclareMathOperator{\vol}{vol}
\DeclareMathOperator{\divv}{div}

\newtheorem{theorem}{Theorem}[section]

\newtheorem{lemma}[theorem]{Lemma}
\newtheorem{corollary}[theorem]{Corollary}

\theoremstyle{definition}

\newtheorem{remark}[theorem]{Remark}
\newtheorem{Example}[theorem]{Example}    
\newtheorem{example}{Example}

\theoremstyle{remark}

\numberwithin{equation}{section}

\begin{document}

\title[Rotationally invariant CMC hypersurfaces in $\hhn$]
{Rotationally invariant hypersurfaces \\ with constant mean curvature in the Heisenberg group
$\hhn$}

\author[M.~Ritor\'e]{Manuel Ritor\'e}
\address{Departamento de Geometr\'{\i}a y Topolog\'{\i}a \\
Universidad de Granada \\ E--18071 Granada \\ Espa\~na}
\email{ritore@ugr.es}

\author[C.~Rosales]{C\'esar Rosales}
\address{Departamento de Geometr\'{\i}a y Topolog\'{\i}a \\
Universidad de Granada \\ E--18071 Granada \\ Espa\~na}
\email{crosales@ugr.es}

\date{April 21, 2005}

\thanks{Both authors have been supported by MCyT-Feder research
project MTM2004-01387.}



\begin{abstract}
In this paper we study sets in the $n$-dimensional Heisenberg group
$\hhn$ which are critical points, under a volume constraint, of the
sub-Riemannian perimeter associated to the distribution of horizontal
vector fields in $\hhn$.  We define a notion of mean curvature for
hypersurfaces and we show that the boundary of a stationary set is a
constant mean curvature (CMC) hypersurface.  Our definition coincides
with previous ones.

Our main result describes which are the CMC hypersurfaces of
revolution in $\hhn$.  The fact that such a hypersurface is invariant
under a compact group of rotations allows us to reduce the CMC partial
differential equation to a system of ordinary differential equations.
The analysis of the solutions leads us to establish a counterpart in
the Heisenberg group of the Delaunay classification of constant mean
curvature hypersurfaces of revolution in the Euclidean space. Hence we
classify the rotationally invariant isoperimetric sets in $\hhn$.
\end{abstract}

\maketitle

\thispagestyle{empty}

\section{Introduction}
\label{sec:intro}
\setcounter{equation}{0}

In the last years the study of variational questions in sub-Riemannian
geometry has experimented an increasing interest.  In particular, the
recent development of a theory of \emph{minimal surfaces} in this
context has contributed to achieve a better understanding of the
geometry of the \emph{Heisenberg group} $\hhn$ endowed with its
\emph{Carnot-Carath\'eodory distance}.

It is well-known that minimal surfaces arise as extremal points of the
area for variations preserving the boundary of the surface.  In this
paper, we are interested in sets contained in the Heisenberg group
which are critical points of the sub-Riemannian perimeter under a
\emph{volume constraint}.  In order to precise the situation we need
to recall some facts about the Heisenberg group, that will be treated 
in more detail in Section~\ref{sec:preliminaries}.

We denote by $\hhn$ the $n$-\emph{dimensional Heisenberg group}, which
we identify with the Lie group $\mathbb{C}^n\times\rr$ where the
product is given by
\[
[z,t]*[z',t']=[z+z',t+t'+\text{Im}\big(\sum_{i=1}^nz_i\overline{z}'_i\big)],
\]
with $z=(z_1,\dots,z_n)$ and $z'=(z'_1,\dots,z'_n)$.  The Lie algebra
of $\hhn$ is generated by $(2n+1)$ left-invariant vector fields
$\{X_k,Y_k,T:k=1,\dots, n\}$ with one non-trivial bracket relation
given by $[X_k,Y_k]=-2T$.  The $(2n)$-dimensional distribution
generated by $\{X_k,Y_k:k=1,\dots,n\}$ is called the \emph{horizontal
distribution} in $\hhn$.  Usually $\hhn$ is endowed with a structure
of sub-Riemannian manifold by considering the Riemannian metric on the
horizontal distribution so that the basis $\{X_k,Y_k:k=1,\ldots,n\}$
is orthonormal.  This metric allows us to measure the length of
horizontal curves and define the \emph{Carnot-Carath\'eodory distance}
between two points as the infimum of length of horizontal curves
joining both points, see \cite{gromov}.  Since $\hhn$ is a group
one can consider its Haar measure, which turns out to coincide with
the Lebesgue measure in $\rrn$.  From the notions of distance and
volume one can also define the Minkowski content in the group, and the
spherical Hausdorff measure, so that different surface measures are
also given on $\hhn$.  As it is shown in \cite{msc} and
\cite{perimeter}, both notions of ``perimeter" coincide for a set with
$C^2$ boundary.

In this paper we will take a slightly different approach to introduce
the notions of volume and perimeter in the Heisenberg group.  We
consider the left-invariant Riemannian metric
$g=\escpr{\cdot\,,\cdot}$ on $\hhn$ so that
$\{X_k,Y_k,T:k=1,\dots,n\}$ is an orthonormal basis at every point.
It is known that the Carnot-Carath\'edory distance can be approximated
in the Gromov-Hausdorff sense by a sequence of dilated Riemannian
metrics associated to $g$, see \cite{gromov2} and \cite{pansu}.  We
define the volume $\vol(\Om)$ of a Borel set $\Om\subeq\hhn$ as the
Riemannian measure of the set.  The perimeter $\pp(\Om)$ is defined in
the sense of De Giorgi by using horizontal vector fields on $\hhn$,
see \eqref{eq:perimeter1} for a precise definition.  This notion of
perimeter coincides with the one introduced in \cite{cng} and
\cite{perimeter}, and it is a sub-Riemannian analogous of the De
Giorgi perimeter in the Riemannian manifold $(\hhn,g)$.

With the notions of volume and perimeter above, we study in
Section~\ref{sec:meancurvature} \emph{stationary sets} of $\hhn$
which are critical points of the perimeter functional for volume
preserving variations.  As in the Riemannian case, one may expect that
for such a set some geometric quantity defined on the boundary remains
constant.  By using the first variation of perimeter in
Lemma~\ref{lem:dp/dt} we will see that the boundary of a $C^2$
stationary set must have \emph{constant mean curvature} (CMC).  The
mean curvature of a hypersurface $\Sg$ is defined in \eqref{eq:mc} as
the Riemannian divergence relative to $\Sg$ of the \emph{horizontal
normal vector} $\nu_H$ to $\Sg$.  We remark that a notion of mean
curvature in $\hh^1$ for graphs over the $xy$-plane was previously
introduced by S.~Pauls \cite{pauls}.  A more general definition has
been proposed by J.-H.~Cheng, J.-F.~Hwang, A.~Malchiodi and P.~Yang
\cite{chmy}, and by N.~Garofalo and S.~Pauls \cite{gp}.  In
Section~\ref{sec:computation} we expose a method to compute the mean
curvature of a hypersurface which in particular shows that our
definition coincides with the previous ones.

The recent study of CMC hypersurfaces in $\hhn$ has mainly focused on
minimal surfaces in $\hh^1$.  In fact, many classical questions of the
theory of minimal surfaces in $\rr^3$, such as the Plateau problem,
the Bernstein problem, or the global behaviour of properly embedded
surfaces, have been treated in $\hh^1$, see \cite{pauls}, \cite{chmy},
\cite{gp} and \cite{ch}.  These works also provide a rich variety of
examples of minimal surfaces in $\hh^1$.  However, in spite of the
last advances, very few is known about non-zero CMC surfaces in
$\hhn$.  In \cite{chmy} some very interesting facts about the CMC
equation in $\hh^1$, such as the uniqueness of solutions for the
Dirichlet problem or the structure of the singular set, are studied.
As to the examples, the only known complete hypersurfaces with
non-zero CMC are the compact spherical ones
(Example~\ref{ex:nonzerocmc}) described in \cite{monti} and
\cite{leomas}.  These hypersurfaces are characterized as the
minimizers of perimeter under a volume constraint in the class of sets
in $\hhn$ bounded by two radial graphs over the hyperplane $\{t=0\}$,
see \cite{dgn}.

The aim of this paper is to study invariant CMC hypersurfaces in
$\hhn$.  Some examples of invariant minimal surfaces in $\hh^1$ were
previously given by S. Pauls \cite{pauls}.  He considered the
solutions of the minimal surface equation for radial graphs over the
$xy$-plane, and discovered a family of complete surfaces of revolution
which are similar to catenoids of $\rr^3$, see
Example~\ref{ex:minimalsurfaces}.  With the same idea he also
introduced some examples of translationally invariant and helicoidal
minimal surfaces.

In Section~\ref{sec:main} of the paper we focus our attention on CMC
hypersurfaces of revolution about the $t$-axis in $\hhn$.  In this
case we can reduce the CMC partial differential equation to a system
of ordinary differential equations (Lemma~\ref{lem:mcequation}).
Then, a detailed analysis of the solutions leads us to our main result
(Theorem~\ref{th:main}) where we prove a Heisenberg analogous of the
classification by C.~Delaunay \cite{delaunay} of constant mean
curvature hypersurfaces of revolution in $\rr^3$, later extended by
W.-Y.~Hsiang \cite{hsiang} to $\rr^n$.  As a consequence, we deduce
that the only compact, CMC hypersurfaces of revolution about the
$t$-axis in $\hhn$ are the spherical ones in
Example~\ref{ex:nonzerocmc}.  Also that complete minimal hypersurfaces
of revolution are hyperplanes orthogonal to the $t$-axis and the
catenoidal type hypersurfaces in Example~\ref{ex:minimalsurfaces}.
Right cylinders, unduloid type hypersurfaces and non-embedded nodoids
provide new examples of complete, periodic hypersurfaces with non-zero
CMC in $\hhn$.

In addition to the geometric interest of this work, we believe that
our results could contribute to study the \emph{isoperimetric problem}
in $\hhn$, which consists of finding sets enclosing a given volume
with the least possible perimeter.  It was proved by G. P.~Leonardi
and S.~Rigot \cite{leorig} that solutions to this problem exist: they
are bounded, connected, and satisfy a certain geometric separation
property (\cite[Theorem 2.11]{leorig}).  Though the
\emph{isoperimetric sets} are unknown, we could expect that they are
rotationally invariant about the $t$-axis, up to a left translation.
In case this was proved, then Theorem~\ref{th:main} would show that
the solutions are congruent with the CMC spheres of revolutions in
Example~\ref{ex:nonzerocmc}.  However, the rotationally symmetry of
isoperimetric sets is still an open question.

Finally, we must remark that the study of constant mean curvature
hypersurfaces in the Riemannian manifold $(\hhn,g)$ is also of great
interest.  The complete classification of rotationally invariant
constant mean curvature hypersurfaces in $(\hh^1,g)$ was established
by P.~Tomter~\cite{tomter}.  The case of $(\hhn,g)$ was treated by
C.~Figueroa, F.~Mercuri and R.~Pedrosa \cite{renato}.  These papers
are based on the so-called \emph{reduction technique}, which consists
of considering how the constant mean curvature equation descends to
the Riemannian quotient of $(\hhn,g)$ by a closed subgroup of
isometries, see \cite[\S 2]{renato} and the references therein.  We
remark that our proof of Theorem~\ref{th:main} is not a generalization
to a sub-Riemannian setting of the above mentioned technique.

\section{Preliminaries}
\label{sec:preliminaries}
\setcounter{equation}{0}

In order to introduce the Heisenberg group we will follow the notation
used by P.~Tomter \cite{tomter}.  The \emph{n-dimensional Heisenberg
group} $\hhn$ is the Lie group $(\rrn,*)$, where we consider the usual
differentiable structure in $\rrn\equiv\mathbb{C}^n\times\rr$, and the
product
\[
[z,t]*[z',t']=[z+z',t+t'+\text{Im}\big(\sum_{i=1}^nz_i\,\overline{z}'_i\big)].
\]
For $p=[z,t]\in\hhn$, the \emph{left translation} associated to $p$ is
the diffeomorphism $L_p(q)=p*q$.  A basis of left-invariant vector
fields is given by
\begin{align*}
X_k&=\frac{\ptl}{\ptl x_k}+y_k\,\frac{\ptl}{\ptl t},\qquad k=1,\dots, n,
\\
Y_k&=\frac{\ptl}{\ptl y_k}-x_k\,\frac{\ptl}{\ptl t},\qquad k=1,\dots, n,
\\
T&=\frac{\ptl}{\ptl t},
\end{align*}
where $(x_k,y_k,t)$ are coordinates in $\rrn$.

The \emph{horizontal distribution} in $\hhn$ is the $(2n)$-dimensional
smooth distribution generated by $\{X_k,Y_k:k=1,\dots,n\}$.  The
notation $U_H$ will represent the projection of a vector $U$ to the
horizontal distribution.  A vector field $U$ is called
\emph{horizontal} if $U=U_H$.  Note that
$[X_k,T]=[Y_k,T]=[X_k,X_j]=[Y_k,Y_j]=0$, while
$[X_k,Y_j]=-2\delta_{kj}\,T$, where $\delta_{kj}$ is the Kronecker
delta.  The last equality and Frobenius Theorem imply that the
horizontal distribution is not integrable.

For a $C^1$ hypersurface $\Sg\sub\hhn$ the \emph{singular set} $\Sg_0$
consists of those points $p\in\Sg$ for which the tangent hyperplane
$T_p\Sg$ coincides with the horizontal distribution.  The set $\Sg_0$
is closed and has empty interior in $\Sg$.  Hence, the \emph{regular
set} $\Sg-\Sg_0$ of $\Sg$ is open and dense in $\Sg$.  For any point
$p\in\Sg-\Sg_0$, the tangent hyperplane meets transversally the
horizontal distribution, and so the intersection is
$(2n-1)$-dimensional.

Consider the Riemannian metric $g=\escpr{\cdot\,,\cdot}$ on $\hhn$ so
that $\{X_k,Y_k,T:k=1\ldots n\}$ is an orthonormal basis of $\rrn$ at
every point.  We denote by $|U|$ the modulus of a vector field $U$.
The volume $\vol(\Om)$ of a Borel set $\Om\subeq\hhn$ is the
Riemannian volume of the metric, which in this case coincides with the
Lebesgue measure in $\rrn$.  The \emph{perimeter} of a Borel set
$\Om\subeq\hhn$ is defined as
\begin{equation}
\label{eq:perimeter1}
\pp(\Om)=\sup\,\left\{\int_\Om\divv(U)\,dv: |U|\leq 1\right\},
\end{equation}
where the supremum is taken over $C^1$ \emph{horizontal vector fields}
with compact support on $\hhn$.  In the definition above, $dv$ and
$\divv(\cdot)$ are the Riemannian volume and divergence of the metric,
respectively.  It is not difficult to see that our notion of perimeter
coincides with the sub-Riemannian perimeter in $\hhn$ given in
\cite{cng} and \cite{perimeter}.  A set $\Om$ is said to be of
\emph{finite perimeter} if $\vol(\Om)$ and $\pp(\Om)$ are finite.  We
refer to the reader to \cite{perimeter} for a detailed development
about perimeter and sets of finite perimeter in $\hhn$.

Let $\Om$ be an open set in $\hhn$ bounded by a $C^2$ embedded
hypersurface $\Sg=\ptl\Om$.  We denote by $N$ the unit normal vector
to $\Sg$ in $(\hhn,g)$ pointing into $\Om$.  By using the Riemannian
divergence theorem we obtain
\begin{equation}
\label{eq:perimeter2} \pp(\Om)=\int_\Sg|N_H|\,da,
\end{equation}
where $da$ is the Riemannian measure on $\Sg$.

We shall denote by $D$ the Levi-Civit\'a connection on $(\hhn,g)$.
The following derivatives can be easily computed
\begin{alignat}{2}
\notag D_{X_k}X_j&=D_{Y_k}Y_j=D_{T}T=0, \\
\label{eq:christoffel}
D_{X_k}Y_j&=-\delta_{kj}\,T, \qquad D_{X_k}T=Y_k, \qquad D_{Y_k}T=-X_k, \\
\notag D_{Y_k}X_j&=\delta_{kj}\,T, \qquad \ \ \,D_{T}X_k=Y_k, \qquad D_{T}Y_k=-X_k.
\end{alignat}
For any vector field $U$ on $\hhn$ we define $G(U)=D_UT$.  It follows
that $G(X_k)=Y_k$, $G(Y_k)=-X_k$ and $G(T)=0$, so that $G$ defines a
linear isometry when restricted to horizontal vector fields.  Note
also that
\begin{equation}
\label{eq:conmute} \escpr{G(U),V}+\escpr{U,G(V)}=0,
\end{equation}
for any pair of vector fields $U$ and $V$.

Let $\Sg$ be a $C^2$ hypersurface in $\hhn$, and $N$ a unit normal
vector to $\Sg$.  We can describe the singular set $\Sg_0\sub\Sg$ in
terms of $N_H$, as the set $\{p\in\Sg:N_H(p)=0\}$.  In the regular
part $\Sg-\Sg_0$, we can define the \emph{horizontal unit normal
vector} $\nu_H$ (\cite{dgn}) by
\begin{equation}
\label{eq:nuh}
\nu_H=\frac{N_H}{|N_H|}.
\end{equation}
Consider the unit vector field $Z$ on $\Sg-\Sg_0$ given by
$Z=G(\nu_H)$.  As $Z$ is horizontal and orthogonal to $\nu_H$, we
conclude that $Z$ is tangent to $\Sg$.

Any isometry of $(\hhn,g)$ leaving invariant the horizontal
distribution preserves the perimeter of sets in $\hhn$.  Examples of
such isometries are left translations, which act transitively on
$\hhn$.

In $\hh^1$, the Euclidean rotation of angle $\theta$ about the
$t$-axis given by
\[
r_{\theta}(x,y,t)=(x\cos\theta-y\sin\theta,x\sin\theta+y\cos\theta,t),
\]
is also such a kind of isometry since it transforms the
orthonormal basis $\{X,Y,T\}$ at the point $p$ into the orthonormal
basis $\{(\cos\theta)\,X+(\sin\theta)\,Y,
(-\sin\theta)\,X+(\cos\theta)\,Y, T\}$ at the point $r_{\theta}(p)$. It
is pointed out in \cite{renato} that not all the rotations about the
$t$ axis are isometries of $(\hhn,g)$ for $n\geq 2$.

\section{Stationary sets and constant mean curvature hypersurfaces in $\hhn$}
\label{sec:meancurvature}
\setcounter{equation}{0}

In this section we study sets of $\hhn$ which are critical points
under a volume constraint of the perimeter functional defined in
\eqref{eq:perimeter1}.

Let $\Om$ be a set of finite perimeter in $\hhn$.  Consider a $C^1$
vector field $U$ with compact support on $\hhn$, and denote by
$\{\varphi_t\}_{t\in\rr}$ the associated group of diffeomorphisms.
Let $\Om_t=\varphi_t(\Om)$.  The fa\-mily $\{\Om_t\}$, for $t$ small,
is the variation of $\Om$ induced by $U$.  Let $V(t)=\vol(\Om_t)$ and
$\pp(t)=\pp(\Om_t)$.  We say that the variation \emph{preserves
volume} if $V(t)$ is constant for $t$ small enough.  We say that $\Om$
is \emph{stationary} if $\pp'(0)=0$ for any volume preserving
variation.  In order to describe analytically stationary sets we shall
compute the first variation formula for volume and perimeter.

Suppose that $\Om$ is bounded by a $C^2$ embedded hypersurface
$\Sg=\ptl\Om$.  As the volume considered is the Riemannian one, it is
well-known (\cite{bdce}) that
\[
V'(0)=\int_\Om\divv(U)\,dv=-\int_\Sg u\,da,
\]
where $u=\escpr{U,N}$ is the component of $U$ with respect to the unit
normal vector $N$ to $\Sg$ pointing into $\Om$.

Now we compute the first variation of perimeter.  We need a previous
lemma.

\begin{lemma}
\label{lem:derivatives} Let $\Sg\subset\hhn$ be a $C^2$ hypersurface
and $N$ a unit normal vector to $\Sg$.  Consider a point
$p\in\Sg-\Sg_0$, the horizontal normal $\nu_H$ at $p$ defined in
\eqref{eq:nuh}, and $Z=G(\nuh)$.  Let $\{Z_1,\dots,Z_{2n-1}\}$ be an
orthonormal family of horizontal, tangent vectors to $\Sg$ at $p$ with
$Z_1=Z$.  Then, for any $u\in T_{p}\hhn$ we have
\begin{align}
\label{eq:uno} D_{u}N_{H}&=(D_{u}N)_{H}-\escpr{N,T}\,G(u)-\escpr{N,G(u)}\,T,
\\
\label{eq:dos} u\,(|N_H|)&=\escpr{D_{u}N, \nuh}-\escpr{N,T}\,\escpr{G(u),\nuh},
\\
\label{eq:tres} D_{u}\nu_H&=|N_H|^{-1}\,\sum_{i=1}^{2n-1}(\escpr{D_uN,Z_i}-\escpr{N,T}\,
\escpr{G(u),Z_i})\,Z_i+\escpr{Z,u}\,T.
\end{align}
\end{lemma}

\begin{proof}
Equalities \eqref{eq:uno} and \eqref{eq:dos} are easily obtained by
using that $N_H=N-\escpr{N,T}\,T$.  Let us prove \eqref{eq:tres}.  As
$|\nu_H|=1$ and $\{\nuh,Z_i,T:i=1,\ldots, 2n-1\}$ is an orthonormal basis of $T_p\hhn$,
we get
\[
D_{u}\nuh=\sum_{i=1}^{2n-1}\escpr{D_{u}\nuh,Z_i}\,Z_i+\escpr{D_{u}\nuh,T}\,T.
\]
Note that $\escpr{D_{u}\nuh,T}=-\escpr{\nuh,G(u)}=\escpr{Z,u}$ by
\eqref{eq:conmute}.  On the other hand, by using \eqref{eq:uno} and
the fact that $Z_i$ is tangent and horizontal, we deduce
\[
\escpr{D_{u}\nuh,Z_i}=|N_H|^{-1}\escpr{D_{u}N_H,Z_i}=
|N_H|^{-1}(\escpr{D_uN,Z_i}-\escpr{N,T}\escpr{G(u),Z_i}).
\]
\end{proof}

For a $C^1$ vector field $U$ on a hypersurface $\Sg$, we denote by
$\divv_\Sg U$ the Riemannian divergence of $U$ relative to $\Sg$,
which is given by
$\divv_{\Sg}U(p):=\sum_{1=1}^{2n}\escpr{D_{e_{i}}U,e_{i}}$ for any
orthonormal basis $\{e_{i}:i=1,\ldots, 2n\}$ of $T_{p}\Sg$.
Now, we can prove

\begin{lemma}
\label{lem:dp/dt} Let $\Om$ be a set of finite perimeter in $\hhn$ such that $\Sg=\ptl\Om$ is a
$C^2$ hypersurface. Suppose that $U$ is a $C^1$ vector field on $\hhn$, whose support in $\Sg$
is disjoint from the singular set $\Sg_0$. Then the first derivative at the origin of the
perimeter functional $\pp(t)$ associated to $U$ is given by
\[
\pp'(0)=\int_\Sg u\,(\divv_\Sg\nu_H)\,da,
\]
where $u=\escpr{U,N}$ is the component of $U$ with respect to the unit
normal vector $N$ to $\Sg$ pointing into $\Om$, and $\nu_H$ is the
horizontal normal vector defined in \eqref{eq:nuh}.
\end{lemma}

\begin{proof}
Call $\{\varphi_t\}_{t\in\rr}$ to the group of diffeomorphisms
associated to $U$ and denote $\Sg_t=\varphi_t(\Sg)$.  Let $da_t$ be
the Riemannian measure on $\Sg_t$.  Consider a $C^1$ vector field $N$
whose restriction to $\Sg_t$ coincides with the unit normal vector
pointing into $\Om_t=\varphi_t(\Om)$.  Denote by $U^\top$ and $U^\bot$
the tangent and the normal part of $U$, respectively.  By using
\eqref{eq:perimeter2} and the coarea formula, we have
\[
\pp(t)=\int_{\Sg_t}|N_H|\,da_t=\int_\Sg (|N_H|\circ\var_t)\,|\text{Jac}\,\var_t|\,da,
\]
where $\text{Jac}\,\var_t$ is the Jacobian determinant of the map
$\var_t:\Sg\to\Sg_t$.  Now, we differentiate with respect to $t$, and
we use the known fact that
$(d/dt)|_{t=0}\,|\text{Jac}\,\var_t|=\divv_\Sg U$, to get
\begin{align*}
\pp'(0)&=\int_{\Sg}\{U(|N_H|)+|N_H|\,\divv_{\Sg}U\}\,da
\\
&=\int_\Sg\{U^\bot(|N_H|)+\divv_\Sg(|N_H|\,U)\}\,da
\\
&=\int_{\Sg}\{\divv_{\Sg}(|N_H|\,U^\top)+
U^\bot(|N_H|)+|N_H|\,\divv_{\Sg}U^\bot\}\,da
\\
&=\int_{\Sg}\{U^\bot(|N_H|)+|N_H|\,\divv_{\Sg}U^\bot\}\,da.
\end{align*}
To obtain the last equality we have used that the integral of the divergence of $|N_H|\,U^\top$
vanishes by virtue of the Riemannian divergence theorem.

On the other hand, we can use \eqref{eq:dos} to obtain
\[
U^{\bot}(|N_H|)=\escpr{D_{U^{\bot}}N,\nu_H} -\escpr{N,T}\,\big<G(U^\bot),\nu_H\big>
=-\escpr{\nabla_\Sg u,\nu_H},
\]
since $G(U^{\bot})$ is orthogonal to $\nu_H$ and
$D_{U^\bot}N=-\nabla_\Sg u$.  Here $\nabla_\Sg u$ represents the
gradient of $u$ relative to $\Sg$.  Then, we get
\begin{align*}
U^\bot(|N_H|)+|N_H|\,\divv_{\Sg}U^\bot&=-(\nu_H)^{\top}(u)+u\,|N_H|\,\divv_{\Sg} N
\\
&=-\divv_\Sg\big(u\,(\nu_H)^\top\big)+u\,\divv_\Sg\big((\nu_H)^\top\big)
+u\,\divv_\Sg(|N_H|\,N)
\\
&=-\divv_\Sg\big(u\,(\nu_H)^\top\big)+u\,\divv_\Sg\nu_H.
\end{align*}
As a consequence, we conclude that
\[
\pp'(0)=-\int_\Sg\divv_\Sg\big(u\,(\nu_H)^\top\big)\,da+\int_\Sg u\,(\divv_{\Sg}\nu_H)\,da,
\]
and the proof follows by using the Riemannian divergence theorem and the fact that $u$ has
compact support disjoint from the singular set $\Sg_0$.
\end{proof}

Let $\Sg$ be a $C^2$ hypersurface in $\hhn$, and $N$ a unit normal
vector field to $\Sg$.  We define the \emph{mean curvature} of $\Sg$
with respect to $N$ by equality
\begin{equation}
\label{eq:mc}
-2n\,H(p)=(\divv_{\Sg}\nu_H)(p),\qquad p\in\Sg-\Sg_0.
\end{equation}
We say that $\Sg$ is of \emph{constant mean curvature} (CMC) if $H$ is
constant on $\Sg-\Sg_0$.  In this case we extend $H$ by its constant
value to the whole $\Sg$.  A \emph{minimal hypersurface} in $\hhn$ is
one for which $H=0$. These definitions have sense even for immersed
hypersurfaces.

The first variation of perimeter can be written in terms of the mean
curvature as
\begin{equation}
\label{eq:firstvariation} \pp'(0)=-2n\int_\Sg H\escpr{U,N}\,da.
\end{equation}

\begin{remark}
The first variation of perimeter and the notion of mean curvature were
first given by S.~Pauls~\cite{pauls} for graphs $\Sg=\{t=f(x,y)\}$ in
$\hh^1$, and later extended by J.-H.~Cheng, J.-F.~Hwang, A.~Malchiodi
and P.~Yang \cite{chmy}, and by N.~Garofalo and S.~Pauls~\cite{gp} to
any $C^2$ surface in $\hh^1$.  The case of $\hhn$ has been recently
treated in \cite{dgn}.  In Section~\ref{sec:computation} we will show
that our definition of mean curvature agrees with the previous ones.
\end{remark}

The first variation of perimeter allows us to prove the following
variational property of stationary sets

\begin{corollary}
Let $\Om$ be a set of finite perimeter in $\hhn$ bounded by a $C^2$
hypersurface $\Sg$.  If $\Om$ is stationary, then $\Sg$ has constant
mean curvature.
\end{corollary}

\begin{proof}
Let $u:\Sg\to\rr$ be a $C^1$ mean zero function with compact support
contained in $\Sg-\Sg_0$.  We can construct, as in \cite[Lemma
2.2]{bdce}, a volume preserving variation given by a vector field $U$
such that $\escpr{U,N}=u$.  As $\Om$ is stationary,
\eqref{eq:firstvariation} gives us
\[
2n\int_{\Sg}Hu\,da=0,
\]
and the proof follows since $u$ is an arbitrary mean zero function.
\end{proof}

\begin{remark}
The mean curvature can be defined locally, off of the singular set, in
an immersed hypersurface $\Sg$ in $\hhn$, and even globally if $\Sg$
is two-sided (there is a well defined unit normal vector $N$).  The
boundary area in $\Sg$ may be defined as $\int_{\Sg}|N_{H}|\,da$ and
it can be easily shown that its first derivative equals
$-\int_{\Sg}2nH\,u\,da$ for any variation of $\Sg$ with initial
velocity vector field $uN$ with compact support.  On the other hand,
even if $\Sg$ does not enclose some given volume, the variation of
volume for a displacement of $\Sg$ in the direction of a normal vector
field $uN$ with compact support can be computed and equals
$-\int_{\Sg}u\,da$.  The interested readers are referred to the paper
by J.~L.~Barbosa, M.~do~Carmo and J.~Eschenburg \cite{bdce} for
details in the Riemannian setting.  Hence one can conclude that
immersed hypersurfaces with constant mean curvature are critical
points of the sub-Riemannian boundary area under volume-preserving
variations.
\end{remark}

\section{Computation of the mean curvature and examples}
\label{sec:computation}
\setcounter{equation}{0}

In this section we describe a method to compute the mean curvature defined in \eqref{eq:mc} of
a hypersurface $\Sg$ in $\hhn$. Then, we will apply it to give an explicit expression for the
mean curvature of a graph $t=f(x,y)$ in $\hh^1$. This will allow us to recall two families of
well-known constant mean curvature hypersurfaces in $\hh^1$. Finally, we will compute the mean
curvature of a hypersurface of revolution in $\hhn$ about the $t$-axis.

Consider a $C^2$ immersion $\phi:B\to\hhn$ defined on a $(2n)$-dimensional Riemannian manifold.
Suppose that $N$ is a unit normal vector field to $\Sg=\phi(B)$ in $(\hhn,g)$. Fix a point
$p\in B$ such that $\phi(p)\in\Sg-\Sg_0$, and consider an orthonormal basis
$\{e_1,\dots,e_{2n}\}$ of $T_p B$. Denote $\ptl_j=e_j (\phi)$. The vectors
$\{\ptl_j:j=1,\dots,2n\}$ form a basis of $T_{\phi(p)}\Sg$. Denote by $\nu_H$ the horizontal
normal vector defined in \eqref{eq:nuh}. A unit, horizontal tangent vector to the non-singular
part $\Sg-\Sg_0$ of $\Sg$ is given by $Z=G(\nu_H)$. Take an orthonormal basis
$\{Z_1,\dots,Z_{2n-1}\}$ of horizontal, tangent vectors to $\Sg$ at $\phi(p)$ with $Z_1=Z$. If
we call $S=\escpr{N,T}\,\nu_H-|N_H|\,T$, then it is clear that $\{Z_i,S\}$ is an orthonormal
basis of $T_{\phi(p)}\Sg$, and so the mean curvature \eqref{eq:mc} of $\Sg$ can be computed
as
\[ -2nH=\sum_{i=1}^{2n-1}\escpr{D_{Z_i}\nu_H,Z_i}+\escpr{D_S\nu_H,S}.
\]
By using the expression for $D_{u}\nu_H$ given in Lemma~\ref{lem:derivatives}, we obtain
\begin{equation}
\label{eq:mcespecial} 2nH=|N_H|^{-1}\,\sum_{i=1}^{2n-1}\,\text{II}(Z_i,Z_i),
\end{equation}
where $\text{II}$ is the second fundamental form of $\Sg$ in $(\hhn,g)$ with respect to the
normal $N$.

From the expression above it is easy to see that for the case $n=1$
\[
D_{Z}Z=2H\nuh,
\]
and we deduce that our definition of mean curvature coincides with the
one given in \cite{chmy}.

Denote by $\text{II}_{ij}=\text{II}(\ptl_i,\ptl_j)=-\escpr{D_{\ptl_i}N,\ptl_j}$. It is clear
that
\begin{equation}
\label{eq:piij} \text{II}_{ij}=\escpr{N,D_{e_i}\ptl_j}.
\end{equation}
On the other hand, if the coordinates of $\ptl_j$ with respect to the
left-invariant basis $\{X_k,Y_k,T\}$ are given by
$(x_{kj},y_{kj},t_j)$, then a straightforward computation by using
\eqref{eq:christoffel} shows that the coordinates of $D_{e_i}\ptl_j$
with respect to $\{X_k,Y_k,T\}$ are
\begin{equation}
\label{eq:dij} \left (e_i(x_{kj})-t_i\,y_{kj}-t_j\,y_{ki},e_i(y_{kj})+t_i\,x_{kj}+t_j\,x_{ki},
e_i (t_j)+\sum_{k=1}^n\,(x_{kj}\,y_{ki}-x_{ki}\,y_{kj})\right).
\end{equation}

The calculation of the coefficients $\text{II}_{ij}$ allows us to
compute $\text{II}(Z_i,Z_i)$ so that we can obtain from
\eqref{eq:mcespecial} an explicit expression for the mean curvature of
$\Sg$ in any particular case.

\begin{example}
Let $\Sg$ be the graph of a function $f\in C^2(B)$ over an open set
$B\subeq\rr^2$.  Then, by using \eqref{eq:dij}, \eqref{eq:piij} and
\eqref{eq:mcespecial} for the immersion $\phi(z)=(z,f(z))$, $z\in B$,
we get
\begin{equation}
\label{eq:graphs} 2H=-\,\frac{(f_y+x)^2\,f_{xx}+(f_x-y)^2\,f_{yy}-2\,(f_x-y)\,
(f_y+x)\,f_{xy}}{\{(f_x-y)^2+(f_y+x)^2\}^{3/2}}.
\end{equation}
The expression above coincides with the divergence in $\rr^2$ of the
horizontal vector field $\nu_H$ projected to $\rr^2$.  It follows that our definition of
mean curvature extends the one given by S. Pauls \cite{pauls}.
\end{example}

Now, we will show some known examples of CMC surfaces in $\hhn$.  Of
course, hyperplanes in $\rrn$ are minimal hypersurfaces in $\hhn$.
The construction of examples has been mainly focused on minimal
surfaces in $\hh^1$, see \cite{pauls}, \cite{ch} and \cite{gp}.  We
are interested in the following ones

\begin{Example}[Catenoidal surfaces in $\hh^1$]
\label{ex:minimalsurfaces} For any $E>0$, let us consider the hyperboloid of revolution $\Sg$
in $\rr^3$ defined in coordinates $[z,t]\in\hh^1$, by equality
\[
t=\pm\sqrt{E^2\,(|z|^2-E^2)},\qquad |z|\geq E.
\]
By using the expression \eqref{eq:graphs} for the mean curvature of a graph, it is easy to
check that $\Sg$ is a smooth, rotationally invariant, minimal surface in $\hh^1$. These
surfaces were characterized by S. Pauls~\cite[Section 4]{pauls} as the unique minimal surfaces
given by the unions of two radial graphs over the $xy$-plane.
\end{Example}

There are no many examples of complete surfaces in $\hhn$ with non-zero constant
mean curvature. The best known are the next ones

\begin{Example}[Spherical hypersurfaces in $\hhn$]
\label{ex:nonzerocmc} For any $H>0$ consider the hypersurface $S_H$ in
$\hhn$ defined in coordinates $[z,t]$, for $z\in\mathbb{C}^n$, by
\[
t=\pm\,\frac{1}{2H^2}\,\{H|z|\,\sqrt{1-H^2\,|z|^2}+\arccos{(H|z|)}\},\qquad |z|\leq\frac{1}{H}.
\]
The hypersurface $S_H$ is compact and homeomorphic to a
$(2n)$-dimensional sphere.  It has two singular points on the
$t$-axis.  In \cite{dgn} it is shown that $S_H$ is $C^2$ but not $C^3$
around the singular points.  It was proved in \cite{leomas} that $S_H$
has constant mean curvature $H$.  These hypersurfaces appeared in
\cite{monti} and \cite{leomas} as the solutions of the restricted
isoperimetric problem in $\hhn$ consisting of finding, for fixed
volume, a minimum of the perimeter functional $\pp(\cdot)$ in the
class of sets in $\hhn$ bounded by two symmetric radial graphs over
the hyperplane $\{t=0\}$.  Recently, D.~Danielli, N.~Garofalo and
D.~Nhieu \cite[Thm.~14.6]{dgn} have proved that they are also solutions of the
isoperimetric problem restricted to a wider class of sets with
cylindrical symmetry.
\end{Example}

Finally, we shall compute the mean curvature of a rotationally
invariant hypersurface in $\hhn$.

Let $\Sg$ be a $C^2$ hypersurface in $\hhn$ which is invariant under
the group of rotations in $\rrn$ about the $t$-axis.  Denote by
$\gamma(s)=(x(s),t(s))$, $s\in I$, the generating curve of $\Sg$ in the
half-plane $\{xt\,(=x_1t):x\geq 0\}$.  We parameterize the
hypersurface $\Sg$ in cylindrical coordinates by the immersion
$\phi:B\to\hhn$ given by $\phi(s,\omega)=(x(s)\,\omega,t(s))$, where
$B$ is the manifold $I\times\mathbb{S}^{2n-1}$ endowed with the
Euclidean metric of $\rr^{2n}$.  We take the unit normal vector to
$\Sg$ in $(\hhn,g)$ whose coordinates with respect to $\{X_k,Y_k,T\}$
are
\begin{equation}
\label{eq:normala} \frac{(xx'\omega_{n+k}-t'\omega_k,-xx'\omega_k-t'\omega_{n+k},x')}
{\sqrt{|\gamma'|^2+x^2(x')^2}},\quad\text{whenever }x>0,
\end{equation}
where $\omega=(\omega_{k},\omega_{n+k})\in\rr^n\times\rr^n$.
In particular, we deduce that the singular set $\Sg_0$ of $\Sg$ is contained on the $t$-axis.
Now, choose a point $p=(s,\omega)\in B$ and a geodesic frame $\{u_2,\dots,u_{2n}\}$ of
$\mathbb{S}^{2n-1}$ around $\omega$, with $u_2(\omega)=(-\omega_{n+k},\omega_k)$. We take the
orthonormal basis $e_1=(1,0)$ and $e_j=(0,u_j)$ of $T_p B=\rr\times T_\omega\mathbb{S}^{2n-1}$.
Note that the coordinates of $\ptl_j=e_j (\phi)$ with respect to the basis $\{X_k,Y_k,T\}$ are
\[
\ptl_1=(x'\omega_k,x'\omega_{n+k},t'),\quad\ptl_2=(-x\omega_{n+k},x\omega_k,x^2),
\quad\ptl_j=(x(u_j)_k,x(u_j)_{n+k},0), \quad j\geq 3.
\]
Then, we can compute from \eqref{eq:piij}, \eqref{eq:dij} and \eqref{eq:normala} the
coefficients $\text{II}_{ij}$ of the second fundamental form of $\Sg$ with respect to $N$,
resulting
\begin{alignat*}{2}
\text{II}_{11}&=\frac{x't''-x''t'-2x(x')^2t'}{\sqrt{|\gamma'|^2+x^2(x')^2}},\qquad
\text{II}_{22}&=&\frac{xt'\,(1+2x^2)}{\sqrt{|\gamma'|^2+x^2(x')^2}},\qquad
\text{II}_{ii}=0,\quad i=3,\dots,2n,
\\
\text{II}_{1j}&=\frac{x(t')^2-x^3(x')^2}{\sqrt{|\gamma'|^2+x^2(x')^2}}\,\,\delta_{1j},\qquad \
\ \ \text{II}_{ij}&=&0, \quad\quad i,j\in\{2,\dots,2n\}, \ i\neq j.
\end{alignat*}
On the other hand, we consider the following orthonormal basis $\{Z_i:i=2,\dots,2n\}$ of
horizontal, tangent vectors to $\Sg$
\[
Z_2=Z=\frac{x\,\ptl_1-(t'/x)\,\ptl_2}{\sqrt{x^2(x')^2+(t')^2}},\qquad
Z_i=\frac{\ptl_i}{x},\quad i=3,\dots,2n.
\]
Finally, it is easy to check from \eqref{eq:mcespecial} that the mean
curvature of $\Sg$ with respect to $N$ in the regular set is the
following
\begin{equation}
\label{eq:primera} 2nH=\frac{x^3\,(x't''-x''t')+
(2n-1)\,(t')^3+2\,(n-1)\,x^2(x')^2t'}{x\,\{x^2(x')^2+(t')^2\}^{3/2}}.
\end{equation}

\begin{example}
When the generating curve $\gamma$ is the graph of a function $x(t)$
defined on the $t$-axis, then equation \eqref{eq:primera} becomes
\[
2nH=\frac{(2n-1)-x^3x''+2\,(n-1)\,x^2(x')^2}{x\,\{1+x^2\,(x')^2\}^{3/2}}.
\]
In particular, a cylinder of radius $r>0$ about the $t$-axis has constant mean curvature
$H=\frac{2n-1}{2nr}$.
\end{example}

\section{Classification of rotationally invariant cmc hypersurfaces in $\hhn$}
\label{sec:main}
\setcounter{equation}{0}

In this section we study in detail the solutions of the CMC equation for hypersurfaces
of revolution in $\hhn$ about the $t$-axis.

Let $\Sg$ be a $C^2$ hypersurface in $\hhn$ which is invariant under the group of rotations in
$\rrn$ about the $t$-axis. Denote by $\gamma$ the generating curve of $\Sg$ in the half-plane
$\{xt \,(=x_1t):x\geq 0\}$. We parameterize the curve $\gamma=(x,t)$ by arc-length $s\in I$. Let $N$ be
the unit normal vector to $\Sg$ given in \eqref{eq:normala}, and $H$ the mean curvature of $\Sg$ with
respect to $N$. Denote by $\sg(s)$ the angle between the tangent vector $\gamma'(s)$ and
$\frac{\ptl}{\ptl t}$. It is clear that $x'=\sin\sg$ and $t'=\cos\sg$. By substituting these equalities
into \eqref{eq:primera}, we deduce the following

\begin{lemma}
\label{lem:mcequation}
The generating curve $\gamma=(x,t)$ of a rotationally invariant hypersurface
in $\hhn$ with constant mean curvature $H$, satisfies the following system of
ordinary differential equations
\[
(*)_H \quad \left\{
\begin{aligned}
x'&=\sin\sg,
\\
t'\,&=\cos\sg,
\\
\sg'&=(2n-1)\,\,\frac{\cos^3\sg}{x^3}+
2\,(n-1)\,\,\frac{\sin^2\sg\,\cos\sg}{x}-2nH\,\,\frac{(x^2\sin^2\sg+\cos^2\sg)^{3/2}}{x^2},
\end{aligned}
\right.
\]

\noindent whenever $x>0$. Moreover, the system above has a first integral:
the function given by

\vspace{-0.3cm}
\begin{equation}
\label{eq:energy} \frac{x^{2n-1}\cos\sg}{\sqrt{x^2\sin^2\sg+\cos^2\sg}}-Hx^{2n}
\end{equation}
is constant along any solution $(x,t,\sg)$.
\end{lemma}

\vspace{0,1cm}\noindent\textbf{Remarks.} 1. Note that the system $(*)_H$ has a singularity for
$x=0$. We will show that the possible contact between a solution $(x,t,\sg)$ and the $t$-axis is
perpendicular. This means that the generated hypersurface $\Sg$ is of class $C^1$ around the
$t$-axis.

2. The existence of a first integral follows from Noether's theorem in the Calculus of
Variations (\cite{gh}) by taking into account that the translations along the $t$-axis preserve the
solutions of $(*)_H$. It can also be obtained by using the arguments given by N.~Korevaar,
R.~Kusner and B.~Solomon, see \cite[p. 480]{KKS}.

3. The constant value $E$ of the function \eqref{eq:energy} will be called the \emph{energy} of
the solution $(x,t,\sg)$. Notice that
\begin{equation}
\label{eq:energia} x^{2n-1}\cos\sg=(E+Hx^{2n})\,\sqrt{x^2\sin^2\sg+\cos^2\sg}.
\end{equation}
The equation above clearly implies
\begin{equation}
\label{eq:cuadrado} (x^{4n-2}-(E+Hx^{2n})^2)\,\cos^2\sg=(E+Hx^{2n})^2\,x^2\sin^2\sg,
\end{equation}
from which we deduce the inequality
\begin{equation}
\label{eq:cota} x^{2n-1}\geq |E+Hx^{2n}|.
\end{equation}

In the following results we gather some elementary properties of the solutions of $(*)_H$.

\begin{lemma}
\label{lem:properties} Let $(x(s),t(s),\sg(s))$ be a solution of $(*)_H$ with energy $E$. Then, we have
\begin{itemize}
\item[(i)] The solution can be translated along the $t$-axis. More precisely,
$(x(s),t(s)+t_0,\sg(s))$ is a solution of $(*)_H$ with energy $E$ for any $t_0\in\rr$.
\item[(ii)] If $x'(s_0)=0$, then the solution is symmetric with respect to the line
$\{t=t(s_0)\}$. As consequence, we can continue a solution by reflecting across the critical
points of $x(s)$.
\item[(iii)] The curve $(x(s_0-s),t(s_0-s),\pi+\sg(s_0-s))$ is a solution of $(*)_{-H}$
with energy $-E$.
\end{itemize}
\end{lemma}

\begin{lemma}
\label{lem:functions} Let $(x(s),t(s),\sg(s))$ be a solution of $(*)_H$. If $\cos\sg(s_0)\neq
0$, then the coordinate $x$ is a function over a small $t$-interval around $t(s_0)$. Moreover
\begin{equation}
\label{eq:dx/dt} \frac{dx}{dt}=\tan\sg,\qquad\frac{d^2x}{dt^2}=\frac{\sg'}{\cos^3\sg},
\end{equation}
where $\sg'$ is the derivative of $\sg$ with respect to $s$.
\end{lemma}

The first integral \eqref{eq:energy} allows us to give the complete description of the
solutions of $(*)_H$. They are of the same types as the ones obtained by C.~Delaunay
\cite{delaunay} when he studied constant mean curvature surfaces of revolution in $\rr^3$.

\begin{figure}[ht]
\centerline{\includegraphics[width=11cm]{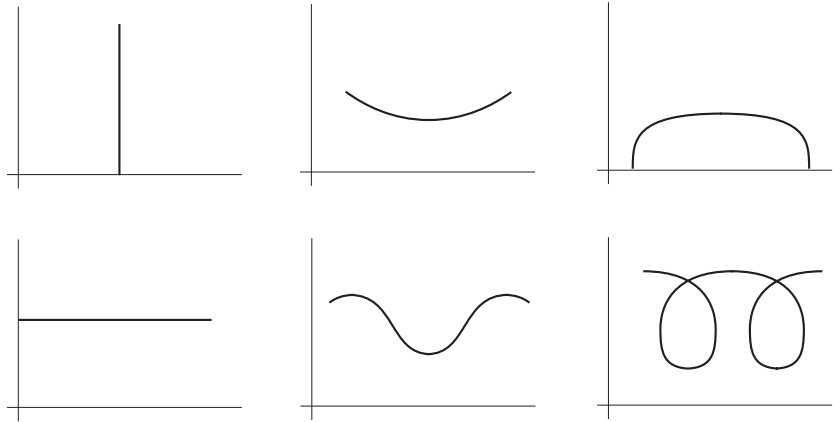}}
\caption{The different types for the generating curve of a
rotationally invariant CMC hypersurface in $\hhn$: hyperplane,
catenoid, sphere, cylinder, unduloid and nodoid.} 
\label{fig:delaunay}
\end{figure}

\begin{theorem}
\label{th:main}
Let $\gamma$ be a complete solution of the system $(*)_H$ with energy $E$.
Then, the generated $\Sg$ is a constant mean curvature hypersurface in $\hhn$ of one of the
following types \emph{(}see Figure \ref{fig:delaunay}\emph{)}
\begin{itemize}
\item[(i)] If $H=0$ and $E=0$ then $\gamma$ is a straight line orthogonal to the $t$-axis and
$\Sg$ is a Euclidean hyperplane.
\item[(ii)] If $H=0$ and $E\neq 0$ we obtain an embedded $\Sg$ of catenoidal type. Moreover,
the resulting hypersurface is contained inside a slab of $\rrn$ when $n\geq 2$.
\item[(iii)] If $H\neq 0$ and $E=0$ then $\Sg$ is a compact hypersurface homeomorphic to a sphere.
\item[(iv)] If $EH>0$ then $\gamma$ is a periodic graph over the $t$-axis. The generated $\Sg$
is a cylinder or an embedded hypersurface of unduloid type.
\item[(v)] If $EH<0$ then $\gamma$ is a locally convex curve and $\Sg$ is a nodoid type hypersurface,
which has selfintersections.
\end{itemize}
\end{theorem}

\begin{proof}
Let $\gamma=(x,t,\sg)$ be a complete solution of $(*)_H$ with energy $E$. By
Lemma~\ref{lem:properties} (i) we can suppose that $\gamma$ is defined over an interval $I$
containing the origin, and that the initial conditions of
$\gamma$ are $(x_0,0,\sg_0)$.  To prove Theorem~\ref{th:main} we
distinguish several cases
depending on the signs of $H$ and $E$. We begin by studying minimal surfaces.

\vspace{0,1cm} \noindent $\bullet$ $H=0$, $E=0$.  From
\eqref{eq:energia} we get $\cos\sg\equiv 0$ and so, by
Lemma~\ref{lem:properties} (iii), we can admit that $\sg\equiv \pi/2$.
It follows that $t\equiv 0$ and $x(s)=s+x_0$.  We conclude that the
solution is a half-line meeting the $t$-axis orthogonally.  The
generated hypersurface $\Sg$ is a hyperplane.

\vspace{0,1cm} \noindent $\bullet$ $H=0$, $E\neq 0$.  We can suppose
that $E>0$ by Lemma~\ref{lem:properties} (iii).  By using
\eqref{eq:cota} we see that $x^{2n-1}\geq E$ and so, the solution does
not approach the $t$-axis.  Moreover, the solution is defined on the
whole real line since $(x',t',\sg')$ is bounded.  By translating the
solution along the $t$-axis we can admit that the initial conditions
of $\gamma$ are $(E^{1/(2n-1)},0,0)$.  From \eqref{eq:energia} we get
$\cos\sg>0$.  It follows from Lemma~\ref{lem:functions} that the $x$-coordinate of
$\gamma$ is a function of $t$ around the origin.  By
$\eqref{eq:dx/dt}$ and $(*)_H$ we have
\[
\frac{d^2x}{dt^2}=\frac{2n-1}{x^3}+\frac{2\,(n-1)}{x}\,\bigg (\frac{dx}{dt}\bigg )^2.
\]
In the case $n=1$ we can integrate this differential equation
obtaining $x(t)=E^{-1}\,\sqrt{t^2+E^4}$, which is a catenoidal type
surface as in Example~\ref{ex:minimalsurfaces}.  Consider the case
$n\geq 2$.  By the symmetry of the solutions of $(*)_H$
(Lemma~\ref{lem:properties} (ii)) we only have to describe the curve
$\gamma(s)$ for $s>0$.  From $(*)_H$ we see that $\sg'>0$, and so
$\sg(s)\in(0,\pi/2)$ for any $s>0$ since $\cos\sg>0$.  On the other
hand, equation \eqref{eq:dx/dt} implies that $x(t)$ is strictly
increasing and strictly convex.  The uniqueness of the solutions of
$(*)_H$ for given initial conditions ensures that any other solution
with $H=0$ and $E>0$ is a translation along the $t$-axis of the graph
of the function $x(t)$ described above.

Let us see that the generated $\Sg$ is contained in a slab of $\rrn$ when $n\geq 2$. Call
$t_\infty=\lim_{s\to +\infty}t(s)$. The fact that $\sin\sg>0$ for $s>0$ implies that the
$t$-coordinate of $\gamma$ is a function of $x\in (E^{1/(2n-1)},+\infty)$. This function
satisfies $dt/dx=\cot\sg$ by \eqref{eq:dx/dt}. An explicit expression for $\cot\sg$ is obtained
by using \eqref{eq:cuadrado}. It follows that
\[
t_\infty=\int_{E^{1/(2n-1)}}^{+\infty}\,\frac{Ex}{\sqrt{x^{4n-2}-E^2}}\,\,dx.
\]
Finally, a comparison with the functions $x^{4n-3}\,(x^{4n-2}-E^2)^{-1/2}$ and $x^{2-2n}$
shows that the integral above converges only if $n\geq 2$.

\vspace{0,1cm} Now, we shall assume that $H\neq 0$. By Lemma~\ref{lem:properties} (iii) we
can suppose $H>0$. We discuss different cases attending to the sign of $E$.

\vspace{0,1cm} \noindent $\bullet$ $H>0$, $E=0$. By equation \eqref{eq:energia} we get
$\cos\sg>0$ on $I$ and so, $\gamma$ is a graph over the $t$-axis. From \eqref{eq:cota} we have
$x\leq 1/H$, so that the solution could approach the $t$-axis. Moreover, by using
\eqref{eq:energia}, \eqref{eq:cuadrado} and $(*)_H$, we obtain
\[
\sg'=\frac{\cos^3\sg}{H^2x^5}\,\,(H^2x^2-2)<0,
\]
which implies that the angle is strictly decreasing along the solution $\gamma$.

As the solutions are invariant by translations along the $t$-axis, we can suppose that the initial
conditions of $\gamma$ are $(1/H,0,0)$. Call $\beta=\sup I$ and
$t_\beta=\lim_{s\to\beta^-}t(s)$. By Lemma~\ref{lem:properties} (ii) we only have to
study the function $x(t)$ for $t\in (0,t_\beta)$. As $\sg'<0$ and $\cos\sg>0$, we deduce that
$\sg\in (-\pi/2,0)$ on the interval $(0,\beta)$. In particular, we get by \eqref{eq:dx/dt} that
$x(t)$ is strictly decreasing and strictly concave on $(0,t_\beta)$. This entails
$t_\beta<+\infty$ since $x(t)$ is bounded. By using the same reasoning with $x(s)$ we prove
$\beta<+\infty$. It follows that $x(s)\to 0$ when $s\to\beta^-$; otherwise, we could continue
the complete solution $\gamma(s)$ for $s>\beta$, a contradiction. By \eqref{eq:energia} we have
that $\gamma$ meets the axis orthogonally. These arguments show that the generated $\Sg$ is
a compact hypersurface of spherical type.

We finally see that $\gamma$ coincides with the generating curve of the sphere $S_H$ given in
Example~\ref{ex:nonzerocmc}. Note that $\sin\sg<0$ on $(0,\beta)$. This implies that the
$t$-coordinate of $\gamma(s)$, for $s\in (0,\beta)$, is a function of $x\in(0,1/H)$. Moreover,
by \eqref{eq:dx/dt} we obtain $dt/dx=\cot\sg<0$. By computing $\cot\sg$ from
\eqref{eq:cuadrado}, we have
\[
\frac{dt}{dx}=\frac{-Hx^2}{\sqrt{1-H^2x^2}},\quad x\in(0,1/H).
\]
We can integrate the equality above to conclude that
\[
t(x)=\frac{1}{2H^2}\,\{Hx\sqrt{1-H^2x^2}+\arccos{(Hx)}\},\quad x\in (0,1/H),
\]
which proves the claim.

\vspace{0,1cm} \noindent $\bullet$ $H>0$, $E>0$. By \eqref{eq:energia} we see that $\cos\sg>0$
on $I$ and so, $\gamma$ is a graph over the $t$-axis. Equation \eqref{eq:cota} gives
$Hx^{2n}-x^{2n-1}+E\leq 0$, which implies that $x(s)$ is bounded and the solution does not meet
the $t$-axis. By $(*)_H$ it follows that $(x',t',\sg')$ is bounded and so, the solution
is defined on the whole real line.  Call $x_1=\inf x(s)$ and $x_2=\sup x(s)$. The values $x_1$
and $x_2$ coincide with the positive zeroes of the polynomial $Hy^{2n}-y^{2n-1}+E$. We consider
two cases

\vspace{0,1cm} \noindent (i) $x_1=x_2$. We obtain $x\equiv r>0$, $\sg\equiv 0$ and $t(s)=s$. The
generated hypersurface is a cylinder of radius $r=(2n-1)/(2nH)$ about the $t$-axis.

\vspace{0,1cm} \noindent (ii) $x_1<x_2$. After a translation along the
$t$-axis, we can suppose that the initial conditions of $\gamma$ are
$(x_1,0,0)$. Call $t_\infty=\lim_{s\to +\infty}t(s)$. By the symmetry of the solutions of
$(*)_H$ it suffices to study the function $x(t)$ for $t\in (0,t_\infty)$. It is clear that
$\sg'(0)\geq 0$; in fact, as we assume $x_1<x_2$, we have that
$\sg'(0)>0$. In particular, by \eqref{eq:dx/dt} we get that $x(t)$ is strictly increasing and
strictly convex on a small interval to the right of the origin. We claim that there exists a
first value $s_1>0$ such that $\sg'(s_1)=0$. Otherwise, we would deduce $\sg(s)\in (0,\pi/2)$
for any $s>0$, which implies by $(*)_H$ that
$x(s)$ is strictly increasing and strictly convex on $(0,+\infty)$, a contradiction since
$x(s)$ is bounded. Call $x_0=x(s_1)$, $t_1=t(s_1)$ and $\sg_1=\sg(s_1)$. The definition of
$s_1$ implies that $\sg\in (0,\pi/2)$ and $\sg'>0$ on $(0,s_1)$. As a consequence

(a) The graph $x(t)$ is strictly increasing and strictly convex on $(0,t_1)$.

On the other hand, by using \eqref{eq:energia}, \eqref{eq:cuadrado} and $(*)_H$, we have
\begin{equation}
\label{eq:sigmaprima1} \sg'=\left[\frac{\cos^3\sg}{x^3\,(E+Hx^{2n})^3}\right]\,p(x),
\end{equation}
where $p$ is the polynomial given by $p(y)=(E+Hy^{2n})^3-2Hy^{6n-2}+2(n-1)Ey^{4n-2}$. As
consequence $p(x_1)>0$ since $\sg'(0)>0$. Moreover, $\sg'(s)=0$ if and only if $p(x(s))=0$. It
is not difficult to see that $p'(y)\neq 0$ for any $y\in [x_1,x_2]$ with $p(y)=0$. Hence,
$p'(x_0)<0$.

It is clear that $x'(s_1)=\sin\sg_1>0$. In particular, $x(s)\neq x_0$ for any $s\neq s_1$
close enough to $s_1$. A straightforward computation shows that
\[
\sg''(s_1)=\left[\frac{\sin\sg_1\,\cos^3\sg_1}{x_0^3\,(E+Hx_0^{2n})^3}\right]\,p'(x_0)<0.
\]
We deduce by using \eqref{eq:dx/dt} that $x(t)$ is strictly concave on a small
interval to the right of $t_1$.

Now, we claim the existence of a first value $s_2>s_1$ such that $\sg(s_2)=0$. Otherwise, $\sg(s)\in
(0,\pi/2)$ for any $s>s_1$, and so $x(s)$ would be strictly increasing and strictly concave on
$(s_1,+\infty)$. Therefore, we would have $x(s)\to x_2$, $x'(s)\to 0$ and
$\sg(s)\to 0$ when $s\to +\infty$. By taking into
account \eqref{eq:sigmaprima1} and the fact that $p(x_2)\neq 0$ we would obtain $\lim_{s\to
+\infty}\sg'(s)\neq 0$, a contradiction. Clearly $x(s_2)=x_2$ by equation
\eqref{eq:energia}; moreover, $\sg'(s_2)<0$ and $p(x_2)<0$. Now, we can prove that $x_0$ is the unique zero of
$p(y)$ in $[x_1,x_2]$. By using Descartes' criterion to count the number of real roots of a
polynomial, we get that $p(y)$ has at most two positive real roots (there are only two
changes of sign in the sequence of coefficients of $p(y)$). But only one positive root is possible
since $p(x_1)>0$,
$p(x_2)<0$, and $p'(y)\neq 0$ for any $y\in [x_1,x_2]$ with $p(y)=0$. Call $t_2=t(s_2)$. By definition of
$s_2$ and the arguments above it follows that $\sg\in (0,\pi/2)$ and $\sg'<0$ on $(s_1,s_2)$. In particular

(b) The graph $x(t)$ is strictly increasing and strictly concave on $(t_1,t_2)$.

As $x'(s_2)=0$ we deduce that the full solution is obtained by successive reflections across
the lines $\{t=k\, t_2\}$ for $k$ an entire. Conclusions (a) and (b) above show that the
generated surface $\Sg$ is similar to an unduloid of $\rrn$.

\vspace{0,1cm} \noindent $\bullet$ $H>0$, $E<0$. In this case we deduce from \eqref{eq:energia}
that the sign of $\cos\sg$ and $E+Hx^{2n}$ is the same. In particular, $\cos\sg=0$ if
and only if $x=x_0=(-E/H)^{1/2n}$. By equation \eqref{eq:cota} we have
$|E+Hx^{2n}|-x^{2n-1}\leq 0$ on $I$. It follows that $x(s)$ is bounded and the solution does
not approach the $t$-axis. Call $x_1=\inf x(s)$ and $x_2=\sup x(s)$. These values are positive zeroes of
the function $|E+Hy^{2n}|-y^{2n-1}$. The case $x_1=x_2$ is not possible; otherwise, the
solution would coincide with a cylinder. By using Descartes' criterion we deduce that the
polynomials $-Hy^{2n}-y^{2n-1}-E$ and $Hy^{2n}-y^{2n-1}+E$ have at most a positive real root.
This implies that $x_0\in (x_1,x_2)$.

On the other hand, by taking into account \eqref{eq:energia} and $(*)_H$, we obtain
\begin{equation}
\label{eq:tela}
\sg'=\frac{\psi}{x^3}\,\,p(x),
\end{equation}
where $p(y)=(E+Hy^{2n})^3-2Hy^{6n-2}+2(n-1)Ey^{4n-2}$, and $\psi$ is the continuous, positive
function defined by
\[
\psi(s)=\left\{
\begin{aligned}
\frac{\cos^3\sg}{(E+Hx^{2n})^3} \ \text{ if }x(s)\neq x_0,
\\
x_0^{6-6n}\qquad\qquad\!\text{ if }x(s)=x_0.
\end{aligned}
\right.
\]
It is easy to check that $p(x_1)<0$ and $p'(y)\leq 0$ for $y\in[x_1,x_2]$. We deduce from \eqref{eq:tela}
that $\sg'<0$ along the solution. After a translation along the $t$-axis we can suppose that the initial
conditions of $\gamma$ are $(x_2,0,0)$. By the symmetry of the
solutions we only have to study the behaviour of $\gamma(s)$ for $s>0$. As $\sg'<0$, we get
$\sg<0$ on $(0,+\infty)$. The fact that $\cos\sg=1$ implies that the $x$-coordinate of $\gamma$ is
a function of $t$ around the origin. By $(*)_H$ and \eqref{eq:dx/dt} it follows that $x(s)$
and $x(t)$ are strictly decreasing and strictly concave on small intervals to the right of the
origin.

Now, the same arguments we used in the case $H>0$, $E>0$ ensure the
existence of a first $s_1>0$ and a first $s_2>s_1$ such that
$\sg(s_1)=-\pi/2$ and $\sg(s_2)=-\pi$.  By \eqref{eq:energia} we see
that $x(s_1)=x_0$ and $x(s_2)=x_1$.  Call $t_1=t(s_1)$ and
$t_2=t(s_2)$.  By the definition of $s_1$ and $s_2$ we have that
$\sg\in (-\pi/2,0)$ on $(0,s_1)$ and $\sg\in (-\pi,-\pi/2)$ on
$(s_1,s_2)$.  As a consequence, the restriction of $\gamma$ to
$[0,s_2]$ consists of two graphs of the function $x(t)$ meeting at
$t=t_1$.  Moreover, by \eqref{eq:dx/dt} we deduce that $x(t)$ is
strictly decreasing and strictly concave on $(0,t_1)$, while it is
strictly decreasing and strictly convex on $(t_2,t_1)$.  As $\{t=0\}$
and $\{t=t_2\}$ are lines of symmetry for $\gamma$, we can reflect
successively to obtain the complete solution, which is periodic.  The
resulting curve is embedded if and only if $t_2=0$; in this case, the
generated $\Sg$ would be compact and homeomorphic to a torus.

Let us see that $\Sg$ has self-intersections. The fact that $\sin\sg<0$ on $(0,s_2)$ implies that we can see
the $t$-coordinate of $\gamma$ as a function of $x\in (x_1,x_2)$. This function satisfies $dt/dx=\cot\sg$ by
\eqref{eq:dx/dt}. An explicit expression for $\cot\sg$ is obtained by using \eqref{eq:cuadrado}. It follows that
\[
t_2=\int_{x_1}^{x_2}\frac{(E+Hx^{2n})\,x}{\sqrt{x^{4n-2}-(E+Hx^{2n})^2}}\,\,dx.
\]
Consider the Riemann surface $R$
associated to the polynomial $w^2=x^{4n-2}-(E+Hx^{2n})^2$.  We may
consider a lift $\alpha$ to $R$ of a Jordan curve $\tilde{\alpha}$ in
the $xt$-plane around the interval $[x_{1},x_{2}]$ so that the only
zeroes of the polynomial $x^{4n-2}-(E+Hx^{2n})^2$ in the interior of
$\tilde{\alpha}$ are $x_{1}$, $x_{2}$. Hence we have
\[
t_{2}=\frac{1}{2}\int_{\alpha}\frac{(E+Hx^{2n})\,x}{w}\,dx.
\]
The function $x^{-2n+2}w$ is holomorphic in a neighbourhood of
$\alpha$. A direct computation yields
\[
d(x^{-2n+2}w)=\bigg(2(n-1)\,x^{-2n+1}(E+Hx^{2n})^2+x^{2n-1}-2nH\,(E+Hx^{2n})\,x
\bigg)\,\frac{dx}{w},
\]
so that
\[
t_{2}=\frac{1}{4nH}\int_{\alpha}\bigg(2(n-1)\,x^{-2n+1}(E+Hx^{2n})^2+x^{2n-1}
\bigg)\,\frac{dx}{w}.
\]
Since the last integrand is an holomorphic one-form in a neighbourhood
of $\alpha$, we finally get
\[
t_{2}=\frac{1}{2nH}\int_{x_{1}}^{x_{2}}\frac{
2(n-1)\,x^{-2n+1}\,(E+Hx^{2n})^2+x^{2n-1}
}{\sqrt{x^{4n-2}-(E+Hx^{2n})^2}}\,\,dx,
\]
which is strictly positive. This allows us to conclude that the generated hypersurface
$\Sg$ is similar to a Euclidean nodoid.
\end{proof}

Theorem~\ref{th:main} provides not only new examples of complete, embedded hypersurfaces in
$\hhn$ with non-zero constant mean curvature, but also the following consequences

\begin{corollary}
\label{cor:minimal} The only minimal hypersurfaces of revolution in $\hhn$ are Euclidean hyperplanes
orthogonal to the $t$-axis and catenoidal type hypersurfaces.
\end{corollary}

\begin{corollary}
\label{cor:notori} The only compact, embedded, rotationally hypersurfaces of constant mean
curvature in $\hhn$ are the spheres $\{S_H\}_{H>0}$ described in Example~\ref{ex:nonzerocmc}.
\end{corollary}

\begin{remark}
Consider the isoperimetric problem in $\hhn$, which consists of finding a minimum for the
perimeter functional $\pp(\cdot)$ in the class of sets in $\hhn$ enclosing a fixed volume. It
was proved by G. Leonardi and S. Rigot \cite{leorig} that the solutions to this problem exist.
Moreover, isoperimetric sets are bounded, connected and satisfy a certain separation property,
see \cite[Theorem 2.11]{leorig}. Though the solutions to this problem are still unknown, we
could expect that they are rotationally invariant about the $t$-axis, up to a left translation.
In case this was proved, then Corollary~\ref{cor:notori} would show that the solutions to the
isoperimetric problem are congruent to the family of spheres $\{S_H\}_{H>0}$ given in
Example~\ref{ex:nonzerocmc}.
\end{remark}

\providecommand{\bysame}{\leavevmode\hbox to3em{\hrulefill}\thinspace}
\providecommand{\MR}{\relax\ifhmode\unskip\space\fi MR }
\providecommand{\MRhref}[2]{%
   \href{http://www.ams.org/mathscinet-getitem?mr=#1}{#2}
}
\providecommand{\href}[2]{#2}

\end{document}